\documentclass[preprint,12pt]{elsarticle}
\usepackage{pgf,tikz}
\usepackage{amsmath,amssymb}
\usepackage{tikz}
\usepackage{color}
\usepackage{xcolor,colortbl}
\usepackage{multirow}
\usepackage{}
\newtheorem{theorem}{Theorem}[section]
\newtheorem{corollary}[theorem]{Corollary}
\newtheorem{conjecture}[theorem]{Conjecture}

\newtheorem{proposition}[theorem]{Proposition}

\newcommand{\pf}{\noindent{\bf Proof.\ }}

\usepackage{subfigure}
\usepackage{float}

\begin{document}

\begin{frontmatter}

\title{On Maximum Induced Forests of the Balanced Bipartite Graphs}

\author[a]{{Ali Ghalavand}\corref{mycorrespondingauthor}}
\cortext[mycorrespondingauthor]{Corresponding author}
\ead{alighalavand@nankai.edu.cn}
\author[a]{{Xueliang Li}}

\address[a]{Center for Combinatorics and LPMC, Nankai University, Tianjin 300071, China}

\begin{abstract}
Let $ \mathcal{B} $ be a balanced bipartite graph with two parts, $ V_1 $ and $ V_2 $, each containing $ n $ vertices, resulting in a total of $ 2n $ vertices. Recently, Wang and Wu conjectured that if the minimum degree of $ \mathcal{B} $, denoted as $ \delta(\mathcal{B}) $, is greater than or equal to $ \frac{n}{2} + 1 $, then the largest order of an induced forest in $ \mathcal{B} $ is equal to $ n + 1 $. In this paper, we prove this conjecture and show that the condition on the minimum degree cannot be relaxed in general terms. Furthermore, we determine that if $ \delta(\mathcal{B}) \geq \frac{n}{2} + 1 $, then any subset $ S $ of vertices in $ \mathcal{B} $ that induces a forest of size $ n + 1 $ will satisfy the conditions  $ \min\{|S \cap V_1|, |S \cap V_2|\} = 1$ when $ n $ is odd, and  $ \min\{|S \cap V_1|, |S \cap V_2|\}\in\{1,2,\frac{n}{2} \}$   when $ n $ is even. Additionally, we identify infinitely many balanced bipartite graphs that meet these conditions.

\end{abstract}


\end{frontmatter}

\section{Introduction}
Let $G$ be a finite and simple graph with vertex set $V(G)$ and edge set $E(G)$.  If $S$ is a subset of $V(G)$, then the induced subgraph of $G$ on $S$ is a graph whose vertex set is $S$ and whose edge set is comprised of all edges in $G$ that connect two vertices in $S$. This induced subgraph is denoted by $G[S]$. For a vertex $v$ in $G$, the degree of $v$ is the number of edges in $G$ that connect to $v$. The degree of $v$ is written as $d_G(v)$. The minimum degree of vertices in $G$ is denoted by $\delta(G)$. Let $n$ be a positive integer. A balanced bipartite graph of order $2n$ is a special type of bipartite graph that has two parts, each containing $n$ vertices, where the total number of vertices is $2n$. We consider the notation $\mathcal{B}$ as a balanced bipartite graph of order $2n$, consisting of two parts $V_1$ and $V_2$.

The forest number of a graph $G$ is the size of the largest subset of vertices of $G$ that form an induced forest. We use $f(G)$ to represent the forest number of graph $G$. A decycling set or a feedback vertex set of a graph is a set of vertices whose removal results in a forest. The smallest possible size of a decycling set of $G$ is represented by $\nabla(G)$. Since $f(G)+\nabla(G)=|V(G)|$, it follows that finding  the decycling number of $G$ is equivalent to determining the largest order of an induced forest, as proposed by Erd\"{o}s et al. in 1986 \cite{Erdos1}.

The problem of destroying all cycles in a graph by deleting a set of vertices was first introduced in combinatorial circuit design in 1974 by Johnson \cite{j1}. Since then, it has found applications in various fields, including deadlock prevention in operating systems by Wang et al. in 1985 \cite{w1} and Silberschatz et al. in 2003 \cite{si1}, the constraint satisfaction problem and Bayesian inference in artificial intelligence by Bar-Yehuda et al. in 1998 \cite{ba1}, monopolies in synchronous distributed systems by Peleg in 1998 \cite{pe1} and 2002 \cite{pe2}, the converter's placement problem in optical networks by Kleinberg and Kumar in 1999 \cite{kl1}, and VLSI chip design by Festa et al. in 2000 \cite{fe1}. The feedback vertex set decision problem involves determining, given a graph $G$ and an integer $k$, whether there is a feedback vertex set of size $k$ in $G$. This problem is known to be NP-complete, even when restricted to planar graphs, bipartite graphs, or perfect graphs, as per Karp 1972 \cite{kar1}.

 Akiyama and Watanabe in 1987 \cite{akiyama1} and, independently, Albertson and Haas in 1998 \cite{albertson1} conjectured that every planar bipartite graph on $n$ vertices contains an induced forest on at least $\frac{5n}{8}$ vertices. Motivated by this conjecture, Alon in 2003 \cite{alon1} studied induced forests in sparse bipartite graphs and showed that every bipartite graph on $n$ vertices with an average degree at most $d\geq1$ contains an induced forest on at least $(\frac{1}{2}+{\rm e}^{-bd^2})n$ vertices, for some absolute positive constant $b$. However, there exist bipartite graphs on $n$ vertices with an average degree at most $d\geq1$ that do not contain an induced forest on at least $(\frac{1}{2}+\frac{1}{{\rm e}^{b'\sqrt{d}}})n$ vertices, where $b'$ is an absolute constant. Conlon et al. in 2014 \cite{conlon1} improved Alon's lower bound to $(\frac{1}{2}+d^{-db})n$ for $d\geq2$. Wang et al. in 2017 \cite{wang2} proved that every simple bipartite planar graph on $n$ vertices contains an induced forest on at least $\lceil \frac{4n+3}{7}\rceil$ vertices. Wang and Wu, in their research on the forest number of graphs and their products in 2023 \cite{wang1}, proposed the following conjecture:
\begin{conjecture}\label{conj1}
Let  $\mathcal{B}$ be the balanced bipartite graph on $2n$ vertices. If $\delta(\mathcal{B})\geq\frac{n}{2}+1$,  then $f(\mathcal{B})=n + 1$.
\end{conjecture}

 In the articles \cite{beineke1,shi1}, upper and lower bounds for the forest number of a graph in terms of its order, size, and maximum degree are provided. Beineke and Vandell \cite{beineke1} studied two families of graph products, namely grids and hypercubes. The forest number of hypercubes has been further researched by Focardi and Luccio \cite{focardi1}. The forest number is a significant parameter of a graph and has been widely studied for planar graphs \cite{dross1,dross2,kelly1,le1,Petrusevski1}, regular graphs \cite{punnim1,ren1}, and subcubic graphs \cite{kelly2}.

 In this paper, we examine the properties of the largest subsets of the vertex set of balanced bipartite graphs that form induced forests and prove the following statements.
 \begin{enumerate}
   \item If $\delta(\mathcal{B})\geq\frac{n}{2}+1$,  then $f(\mathcal{B})=n + 1$.
   \item If $n\geq2$, then there is a  balanced bipartite graph $\mathcal{B}$ with $\delta(\mathcal{B})=\lceil \frac{n}{2}\rceil $ and $f(\mathcal{B})=n+2$.
  \item If  $\delta(\mathcal{B})\geq\frac{n}{2}+1$, and $S$ is a subset of the vertex set of $\mathcal{B}$ with cardinality $n+1$ such that $\mathcal{B}[S]$ is a forest, then either $|S\cap V_1|$ or $|S\cap V_2|$ lies within the set $\{1,2,\frac{n}{2}\}$.
  \item (iv) If $n$ is an odd number,  $\delta(\mathcal{B})\geq\frac{n}{2}+1$,  and $S$ is a subset of the vertex set of $\mathcal{B}$ with $|S|=n+1$, then $\mathcal{B}[S]$ is a forest if and only if $\min\{|S\cap\,V_1|,|S\cap\,V_2|\}=1$.
   \item If $\lambda\in\{1,2,\frac{n}{2}\}$, then there are infinitely many balanced bipartite graphs $\mathcal{B}$ with  $\delta(\mathcal{B})\geq\frac{n}{2}+1$, and for a subset $S$ of the vertex set of $\mathcal{B}$, $|S|=n+1$, $\mathcal{B}[S]$ is a forest, and $\min\{|S\cap\,V_1|,|S\cap\,V_2|\}=\lambda$.
   \item If $n$ is an odd number, $\delta(\mathcal{B})\geq\frac{n+1}{2}$, and $\max\{|\{v:v\in\,V_1~\wedge~d_\mathcal{B}(v)=\frac{n+1}{2}\}|,|\{u:u\in\,V_2~\wedge~d_\mathcal{B}(u)=\frac{n+1}{2}\}|\}\leq1$, then $f(\mathcal{B})=n+1$.
   \item If $ k \geq 2 $ is an integer, then there are infinitely many balanced bipartite graphs $ \mathcal{B} $ with $ \delta(\mathcal{B}) = k $ and $ f(\mathcal{B}) = n + 1 $. Additionally, there are infinitely many balanced bipartite graphs $ \mathcal{B} $ with $\delta(\mathcal{B}) = k $ and $f(\mathcal{B}) = n + 2$.
   
 \end{enumerate}

\section{Main Results}

In this section, we present our main findings concerning the structure and size of the largest subsets of vertices in balanced bipartite graphs that form induced forests. The following theorem validates Conjecture \ref{conj1}.

\begin{theorem}\label{th1}
Let $ \mathcal{B} $ be the  balanced bipartite graph on $2n$ vertices. If $\delta(\mathcal{B})\geq\frac{n}{2}+1$, then $f(\mathcal{B})=n + 1$.
\end{theorem}
\pf
Suppose $\mathcal{B}$ is the balanced bipartite graph with $2n$ vertices and parts $V_1$ and $V_2$. Let $v_1 \in V_1$ and $v_2 \in V_2$. Since $\mathcal{B}$ is a bipartite graph with bipartition $(V_1,V_2)$, both $\mathcal{B}[V_1\cup\{v_2\}]$ and $\mathcal{B}[\{v_1\}\cup V_2]$ are forests. Therefore, $f(\mathcal{B})\geq n+1$. Now, suppose $S$ is a subset of the vertex set of $\mathcal{B}$ such that $|S|=n+2$. Assume $S\cap\,V_1=L_1$, $S\cap\,V_2=L_2$, and
$|L_1|\geq|L_2|$. Since $|S|\geq\,n+2$ and $|V_1|=|V_2|=n$, these follow that $|L_1|\geq|L_2|\geq2$. Thus,  one can observe that the set of possible cardinalities of $L_2$ is $\{2, 3, \ldots, \frac{n+1}{2}\}$ when $n$ is odd and $\{2, 3, \ldots, \frac{n+2}{2}\}$ when $n$ is even. Plus, the cardinality of $L_1$ is $n+2-|L_2|$. By using our assumptions we have $\delta(\mathcal{B}) \geq \frac{n}{2} + 1$. Therefore, for every $k$ in $\{2, 3, ..., \frac{n+1}{2}\}$, if $|L_2|=k$, then $|N_\mathcal{B}(x) \cap L_1| \geq\frac{n}{2} + 1-(n-(n+2-k)) =\frac{n}{2} + 3 - k$ for all $x\in L_2$. Thus, for every $k$ in $\{2, 3, ..., \frac{n+2}{2}\}$, if $|L_2|=k$, then
\begin{equation}\label{eq1}
\sum_{x\in\,L_2}|N_\mathcal{B}(x) \cap L_1| \geq\,k(\frac{n}{2}+3-k).
\end{equation}
One can see that the function $g(k)=k(\frac{n}{2}+3-k)$ is continuous on both intervals $[2, \frac{n+1}{2}]$ and $[2, \frac{n+2}{2}]$. Additionally, $\frac{\partial}{\partial\,k}(g(k))=\frac{n}{2} + 3 - 2k$. Therefore, it is well-known that the function $g$ attains its minimum at one of the points $2$, $\frac{n+1}{2}$, and $\frac{n+6}{4}$ in the domain $[2, \frac{n+1}{2}]$, and it attains its minimum at one of the points $2$, $\frac{n+2}{2}$, and $\frac{n+6}{4}$ in the domain $[2, \frac{n+2}{2}]$. By utilizing these statements along with Equation \eqref{eq1}, we can derive the following results:
\[\sum_{x\in\,L_2}|N_\mathcal{B}(x) \cap L_1| \geq\min\{g(2),g(\frac{n+1}{2}),g(\frac{n+2}{2}),g(\frac{n+6}{4})\}\geq\,n+2.\]
This implies that the subgraph $ \mathcal{B}[S] $ contains at least $ n + 2$ edges. So, by a well-known theorem in graph theory, since $\mathcal{B}[S]$ is an induced subgraph of $\mathcal{B}$ by $n+2$ vertices and greater or equal to $n+2$ edges, it is not an acyclic graph. Therefore, the forest number of graph $\mathcal{B}$, $f(\mathcal{B})$, cannot be greater than or equal to $n+2$. However, we have already proved that $f(\mathcal{B})$ is greater than or equal to $n+1$. Thus, $f(\mathcal{B})$ must be equal to $n+1$, as wanted.
\hfill $\Box$\\

The following proposition demonstrates that the minimum degree condition in Theorem \ref{th1} is sharp.

\begin{proposition}\label{pro1}
If $n\geq2$,  then there is a  balanced bipartite graph $\mathcal{B}$ with $2n$ vertices, minimum degree  $\lceil \frac{n}{2}\rceil $, and $f(\mathcal{B})=n+2$.
\end{proposition}
\pf
Let $n$ be an integer greater than or equal to $2$. Consider a complete balanced bipartite graph $K_{n,n}$ with $2n$ vertices divided into two parts, $A$ and $B$. Assume that $\{a_1,a_2\}\subseteq A$ and $B=\{b_1,b_2,\ldots,b_n\}$. Now, consider the balanced bipartite graph $\mathcal{B}$, which is created by removing the edges between $a_1$ and the first $\lfloor \frac{n}{2}\rfloor$ vertices of $B$, as well as the edges between $a_2$ and the last $n-\lceil \frac{n}{2}\rceil$ vertices of $B$. The structure of $\mathcal{B}$ is such that $d_\mathcal{B}(a_1)=d_\mathcal{B}(a_2)=\lceil \frac{n}{2}\rceil$, $d_\mathcal{B}(x)\geq\,n-1$ for $x\in\,V(\mathcal{B})\backslash\,\{a_1,a_2\}$, and $\mathcal{B}[B\cup\,\{a_1,a_2\}]$ is a forest. It follows that $f(\mathcal{B})\geq\,n+2$. However, one can observe that  $f(\mathcal{B})\leq\,n+2$, because of the structure of $\mathcal{B}$. Thus, we conclude that $f(\mathcal{B})=n+2$, as desired.
\hfill $\Box$\\

The following three results demonstrate the structures of the vertex subsets of the balanced bipartite graphs in Theorem \ref{th1} that form maximum forests.

\begin{theorem}\label{th2}
Let $\mathcal{B}$ be a balanced bipartite graph with $2n$ vertices, partitioned into $V_1$ and $V_2$, and $\delta(\mathcal{B})\geq\frac{n}{2}+1$. If $S$ is a subset of the vertex set of $\mathcal{B}$ with $n+1$ vertices such that $\mathcal{B}[S]$ is a forest, then either $|S\cap V_1|$ or $|S\cap V_2|$ lies within the set $\{1,2,\frac{n}{2}\}$.
\end{theorem}
\pf
Assume  $L_1=S\cap V_1$, $L_2=S\cap V_2$, and $|L_1|\geq|L_2|$. We consider the three cases as follows:
\begin{enumerate}
  \item $|L_2|\in\{1,2,\frac{n}{2}\}$. In this case  there is nothing to prove.
  \item $|L_1|=|L_2|=\frac{n+1}{2}$. In this case by some simple calculations,  $|V_1\backslash\,L_1|=|V_2\backslash\,L_2|=\frac{n-1}{2}$. For every vertex $v$ in $V(\mathcal{B})$, it is assumed that $d_\mathcal{B}(v) \geq \frac{n}{2}+ 1$. Therefore, for every vertex $x$ in $L_1\cup\,L_2$, $d_{\mathcal{B}[S]}(x)\geq\frac{n}{2}+1-\frac{n-1}{2}=\frac{3}{2}$. However, for every vertex $x$ in $L_1\cup\,L_2$, $d_{\mathcal{B}[S]}(x)$ is an integer. Thus, for every vertex $x$ in $L_1\cup\,L_2$, $d_{\mathcal{B}[S]}(x)\geq2$ and it follows that $G[S]$ cannot be a forest.
  \item Otherwise, when $n$ is odd, the set of possible cardinalities of $L_2$ is $\{3,4, \ldots, \frac{n-1}{2}\}$ and when $n$ is even, it is $\{3,4, \ldots, \frac{n-2}{2}\}$. The cardinality of set $L_1$ can be calculated as $n+1-|L_2|$. For every vertex $v$ in $V(\mathcal{B})$, it is assumed that $d_\mathcal{B}(v) \geq \frac{n}{2}+ 1$. Therefore, for every $k$ in $\{2, 3, ..., \frac{n-1}{2}\}$, if $|L_2|=k$, then $|N_\mathcal{B}(x) \cap L_1| \geq \frac{n}{2} + 2 - k$ for all $x\in L_2$. Thus, for every $k$ in $\{2, 3, ..., \frac{n-1}{2}\}$, if $|L_2|=k$, then
\[\sum_{x\in\,L_2}|N_\mathcal{B}(x) \cap L_1| \geq \sum_{x\in\,L_2}(\frac{n}{2} + 2 - k)=k(\frac{n}{2}+2-k).\]
 The function $h(k)=k(\frac{n}{2}+2-k)$ is continuous on both intervals $[2, \frac{n-1}{2}]$ and $[2, \frac{n-2}{2}]$. Plus, $\frac{\partial}{\partial\,k}(h(k))=\frac{n}{2} + 2 - 2k$. Therefore, by using a similar approach as in the proof of Theorem \ref{th1}, we obtain the following results:
\[\sum_{x\in\,L_2}|N_\mathcal{B}(x) \cap L_1| \geq\min\{h(3),h(\frac{n-2}{2}),h(\frac{n-1}{2}),h(\frac{n+4}{4})\}\geq\,n+1.\]
 Thus, $\mathcal{B}[S]$ is an induced subgraph of $\mathcal{B}$ with $n+1$ vertices and at least $n+1$ edges. It follows that $\mathcal{B}[S]$ cannot be a forest. 
\end{enumerate}
By combining the three cases mentioned above, we arrive at the result we were aiming for.
\hfill $\Box$\\

The following theorem demonstrates that if subset $S$ in Theorem \ref{th2} satisfies either $|S \cap V_1| = 2$ or $|S \cap V_2| = 2$, then $n$ must be even. 

\begin{theorem}\label{th4}
Let $\mathcal{B}$ and $S$ be the graph and set considered in the last theorem. If  $\min\{|S \cap V_1|, |S \cap V_2|\} = 2$, then  $n$ must be even.
\end{theorem}
\pf
Suppose $S$ is a subset of the vertices of $\mathcal{B}$ that satisfies the three conditions of the theorem. Without loss of generality, let $|S\cap V_1|=2$, $L_1=S\cap V_1=\{u,v\}$, and $L_2=S\cap V_2$. Since $\delta(\mathcal{B})\geq\frac{n}{2}+1$, it follows that both $d_\mathcal{B}(u)$ and $d_\mathcal{B}(v)$ are  greater than or equal to $\frac{n}{2}+1$. Therefore, $|N_\mathcal{B}(u)\cap L_2|$ and $|N_\mathcal{B}(v)\cap L_2|$ are both greater than or equal to $\frac{n}{2}$, since $|L_2|=n-1$. Now, assume that $n$ is odd. Since $|N_\mathcal{B}(u)\cap L_2|$ and $|N_\mathcal{B}(v)\cap L_2|$ are integers, they must both be greater than or equal to $\frac{n+1}{2}$. It follows that $\mathcal{B}[S]$ is an induced subgraph of $\mathcal{B}$ with $n+1$ vertices and at least $n+1$ edges. Therefore, $\mathcal{B}[S]$ is not a forest, which contradicts our assumption that $\mathcal{B}[S]$ is a forest. Thus, we have shown that $n$ must be even, as desired.
\hfill $\Box$\\

The following corollary shows that if $n$ is odd, then subset $S$ in Theorem \ref{th2} induces a forest if and only if $\min\{|S\cap V_1|, |S\cap V_2|\} = 1$.

\begin{corollary}\label{co1}
Let $\mathcal{B}$ and $S$ be the graph and set considered in the last theorem. If $n$ is odd, then $G[S]$ is a forest if and only if $\min\{|S\cap\,V_1|,|S\cap\,V_2|\}=1$.
\end{corollary}
\pf Assume that $ n $ is odd. Since $ \min\{|S \cap V_1|, |S \cap V_2|\} $ is an integer, Theorems \ref{th2} and \ref{th4} imply that $ \min\{|S \cap V_1|, |S \cap V_2|\} = 1 $. Furthermore, for any $ v_1 \in V_1 $ and $ v_2 \in V_2 $, both $ \mathcal{B}[\{v_1\} \cup V_2] $ and $ \mathcal{B}[\{v_2\} \cup V_1] $ are forests. Thus, the condition $ \min\{|S \cap V_1|, |S \cap V_2|\} = 1 $ ensures that $ G[S] $ is a forest, as desired. \hfill $\Box$\\

In Theorem \ref{th2}, we established that the minimum cardinality of the intersection between the subset $ S $ and the partitions $ V_1 $ and $ V_2 $ of the balanced bipartite graph $ \mathcal{B} $ can only be $ 1 $, $ 2 $, or $ \frac{n}{2} $. The following theorem demonstrates that there are infinitely many balanced bipartite graphs $ \mathcal{B} $ that meet the conditions outlined in Theorem \ref{th2}. Additionally, for each case, the minimum cardinality of the intersection of the subset $ S $ with the partitions $ V_1 $ and $ V_2 $ equals to $\lambda$, where $\lambda\in\{1,2,\frac{|V(\mathcal{B})|}{4}\}$.
\begin{theorem}\label{th3}
There are infinitely many balanced bipartite graphs $\mathcal{B}$ that satisfy the following conditions:
\begin{enumerate}
  \item the graph $\mathcal{B}$ has minimum degree of at least $\frac{|V(\mathcal{B})|}{4}+1$;
  \item there is at least one subset $S$ of $V(\mathcal{B})$ with $\frac{|V(\mathcal{B})|}{2}+1$ elements that induces a forest, and its intersection  with at least one of the partitions of $\mathcal{B}$ has cardinality $\lambda$, where $\lambda\in\{1,2,\frac{|V(\mathcal{B})|}{4}\}$.
\end{enumerate}
\end{theorem}
\pf
Suppose $n$ is a positive integer. We can classify the three cases based on the value of $\lambda$ as follows:
\begin{enumerate}
\item $\lambda=1$.  This case holds for any bipartite graph, as proven in Theorem \ref{th1}.
\item $\lambda=2$. In this case, for any positive even number $n$, we can form a bipartite graph $X$ with bipartition $(A, B)$ where $A = \{a_1, a_2\}$, $|B| = n-1$, $d_X(a_1)=d_X(a_2)=\frac{n-2}{2}+1$, and $|N_X(a_1)\cap\,N_X(a_2)|=1$. We also create another bipartite graph $Y$ with bipartition $(F, H)$, where $|F| = n-2$, $H = \{h\}$, $d_Y(f)\geq0$ for all $f\in\,F$, and $d_Y(h)\geq\frac{n}{2}-1$. Finally, we define $\mathcal{B}$ as the bipartite graph with bipartition $(A\cup\,F, B\cup\,H)$ and edge set $E(X)\cup\,E(Y)\cup\{uv:(u\in\,A~\wedge~v\in\,H)~\vee~(u\in\,B~\wedge~v\in\,F)\}$. One can check that $\mathcal{B}$ is a balanced bipartite graph on $2n$ vertices with a minimum degree of at least $\frac{n}{2}+1$, and for $S=A\cup\,B$, the three statements $|S|=n+1$, $\mathcal{B}[S]$ is a forest, and $|S\cap(A\cup\,F)|=2$ hold.
\item $\lambda=\frac{n}{2}$. In this case, for any positive even number $n$, we can form a bipartite graph $X$ with parts $A$ and $B$ such that $|A|=|B|-1=\frac{n}{2}$, $X$ is isomorphic to $P_{n+1}$, and $d_X(a)=2$ for all $a\in\,A$. We also create another bipartite graph $Y$ with parts $F$ and $H$ such that $|F|=|H|+1=\frac{n}{2}$, $d_Y(f)\geq0$ for all $f\in\,F$, and $d_Y(h)\geq1$ for all $h\in\,H$. Finally, we define $\mathcal{B}$ as a bipartite graph with parts $A\cup\,F$ and $B\cup\,H$, and edge set $E(X)\cup\,E(Y)\cup\{uv:(u\in\,A~\wedge~v\in\,H)~\vee~(u\in\,B~\wedge~v\in\,F)\}$. One can see that $\mathcal{B}$ is a balanced bipartite graph on $2n$ vertices with a minimum degree of at least $\frac{n}{2}+1$, and for $S=A\cup\,B$, the three statements $|S|=n+1$, $\mathcal{B}[S]$ is a forest, and $|S\cap(A\cup\,F)|=\frac{n}{2}$ hold.  
\end{enumerate}
By combining the three cases discussed above, we achieve the desired result.\hfill $\Box$\\

In the following theorem, by employing Theorem \ref{th1}, Corollary \ref{co1}, and additional arguments, we show that for any integers $k \geq 2$ and $l \in \{1, 2\}$, there are infinitely many balanced bipartite graphs $G$ with $\delta(G) = k$ and $f(G) = \frac{|V(G)|}{2} + l$.

\begin{theorem}\label{thh1}
Let $k\geq2$ be an integer, and let  $l \in \{1, 2\}$. Then there are infinitely many balanced bipartite graphs $G$ with $\delta(G) = k$ and $f(G) = \frac{|V(G)|}{2} + l$.
\end{theorem}

\pf
Let $n$, $k$, and $l$ be three positive integers, where $n$ is odd and satisfies the conditions $n \geq k - 1 \geq 1$. Additionally, let $l$ take on one of the values in the set $\{1, 2\}$. We can classify the two cases based on the value of $l$ as follows:
\begin{enumerate}
  \item $l=1$. In this case, suppose $H$ is a balanced bipartite graph with $2n$ vertices and parts $V_1$ and $V_2$. Assume that $\delta(H)\geq\frac{n}{2}+1$ and $\max_{v\in V_2}d_H(v)=n$. Let $w\in V_2$ with $d_H(w)=n$ and $\{w_1,w_2,\ldots,w_{k-2}\}\subseteq V_2\backslash \{w\}$. Also, let $P_2$ be a path with two vertices $x$ and $y$. Consider the balanced bipartite graph $G$ with parts $V_1\cup\{x\}$ and $V_2\cup\{y\}$ and edge set $E(H)\cup\{xw_i:i\in\{1,2,\ldots,k-2\}\}\cup\{yu:u\in V_1\}\cup\{xw\}$. One can verify that $\delta(G)=d_G(x)=k$, $|V(G)|=2n+2$, and $f(G)\geq n+2$. However, by Theorem \ref{th1}, $f(H)=n+1$. Therefore, by Corollary \ref{co1}, if $S\subseteq V(G)$, $|S|=n+3$, and $G[S]$ is a forest, then $\min\{|S\cap V_1|,|S\cap V_2|\}=1$ and $\{x,y\}\subseteq S$. We now consider two cases: (i) $S\cap V_1=\{a\}$ and $S\cap V_2=V_2$. In this case, $N_G(y)\cap N_G(w)=\{x,a\}$, which is a contradiction since $G[S]$ is a forest. (ii) $S\cap V_1=V_1$ and $S\cap V_2=\{b\}$. In this case, $|V(G[S])|=n+3$ and $|E(G[S])|\geq n+3$ (because $d_G(y)=n+1$ and $d_G(b)\geq 2$), which is again a contradiction since $G[S]$ is a forest. Thus, $f(G)=n+2=\frac{|V(G)|}{2}+1$.
  \item $l=2$. Let $H$ be the balanced bipartite graph defined in the previous case, with $\delta(H)=\frac{n}{2}+1$. Without loss of generality, let $v\in V_1$, $d_H(v)=\frac{n}{2}+1$, $V_2\backslash N_H(v)=\{w_1,w_2,\ldots,w_{\frac{n}{2}-1}\}$, $U\subseteq V_1$, and $|U|\geq k$. Let $P_2$ be a path with two vertices $x$ and $y$. Now, consider the balanced bipartite graph $G$ with parts $V_1\cup\{x\}$ and $V_2\cup\{y\}$ and edge set $E(H)\cup\{xw_i:i\in\{1,2,\ldots,k-1\}\}\cup\{yu:u\in U\}$. It can be shown that $\delta(G)=d_G(x)=k$ and $|V(G)|=2n+2$. However, by Theorem \ref{th1}, $f(H)=n+1$. Therefore, for any subset of $V(G)$, such as $S$, that makes $G[S]$ a forest, we have $|S|\leq n+3$, because $|S\cap V(H)|\leq n+1$. Suppose $S=V_2\cup\{x,v\}$. In this case, by the structure of $G$, it can be observed that $|N_G(v)\cap N_G(x)|\leq1$. Thus, $G[S]$ is a forest, and so $f(G)\geq n+3$. Since we have already proven that $f(G)\leq n+3$, $f(G)=n+3=\frac{|V(G)|}{2}+2$.
\end{enumerate}
By combining the two cases mentioned above, we arrive at the result we were aiming for.\hfill $\Box$\\

In the following theorem, we relax the minimum degree condition in Theorem \ref{th1} for a class of balanced bipartite graphs.

\begin{theorem}\label{th5}
Let $\mathcal{B}$ be a balanced bipartite graph with $2n$ vertices, $\delta(\mathcal{B})\geq\frac{n+1}{2}$, and parts $V_1$ and $V_2$. If $n$ is odd, and $\max\{|\{v:v\in\,V_1~\wedge~d_\mathcal{B}(v)=\frac{n+1}{2}\}|,|\{u:u\in\,V_2~\wedge~d_\mathcal{B}(u)=\frac{n+1}{2}\}|\}\leq1$, then $f(\mathcal{B})=n+1$.
\end{theorem}
\pf
Suppose $\mathcal{B}$ is a balanced bipartite graph with two parts, $V_1$ and $V_2$. It is straightforward to see that $f(\mathcal{B}) \geq n + 1$. Now, let us consider a subset $S$ of the vertex set of $\mathcal{B}$ such that $|S| = n + 2$, with $S \cap V_1 = L_1$ and $S \cap V_2 = L_2$. We have the condition $|L_1| \geq |L_2|$.
As shown in the proof of Theorem \ref{th1}, the possible cardinalities of $L_2$ are $\{2, 3, \ldots, \frac{n + 1}{2}\}$, and the cardinality of $L_1$ is given by $n + 2 - |L_2|$. It is assumed that $d_\mathcal{B}(v) \geq \frac{n + 1}{2}$ for every vertex $v \in V(\mathcal{B})$ and that $\max\{|\{v:v\in V_1~\wedge~d_\mathcal{B}(v)=\frac{n+1}{2}\}|,|\{u:u\in V_2~\wedge~d_\mathcal{B}(u)=\frac{n+1}{2}\}|\}\leq1$. Therefore, for every $k$ in $\{2, 3, \ldots, \frac{n + 1}{2}\}$, if $|L_2| = k$, then 
\[\sum_{x\in L_2}|N_\mathcal{B}(x) \cap L_1| \geq\, \frac{n+1}{2}+2-k+(k-1)(\frac{n+1}{2}+3-k).\]
By a method similar to that used in part of the proof of Theorem \ref{th1}, we can conclude that 
\[\sum_{x\in L_2}|N_\mathcal{B}(x) \cap L_1| \geq\min\{g(2),g(\frac{n+1}{2}),g(\frac{n+7}{4})\}\geq\,n+2.\] 
 This implies that $\mathcal{B}[S]$ is not an acyclic graph. Consequently, the forest number of graph $\mathcal{B}$ cannot be greater than or equal to $n + 2$. However, we have already established that $f(\mathcal{B}) \geq n + 1$. Therefore, we can conclude that $f(\mathcal{B})$ must be equal to $n + 1$, as desired.
\hfill $\Box$

\section{Concluding Remarks and Future Work}
Let $ \mathcal{B} $ be a balanced bipartite graph consisting of two parts, $ V_1 $ and $ V_2 $, each containing $ n $ vertices, resulting in a total of $ 2n $ vertices with forest number $ f(\mathcal{B})$ and decycling number $ \nabla(\mathcal{B}) $.  This paper investigates the properties of the largest subsets of the vertex set of the balanced bipartite graph $ \mathcal{B} $ that form induced forests. We proved a conjecture by Wang and Wu, which states that if $ \delta(\mathcal{B})\geq\frac{n}{2} + 1 $, then $ f(\mathcal{B}) = n + 1 $. We also show that this condition on the minimum degree cannot be generally relaxed. Furthermore, we determined that if $ \delta(\mathcal{B}) \geq \frac{n}{2} + 1 $, then any subset $ S $ of vertices in $ \mathcal{B} $ that induces a forest of size $ n + 1 $ will satisfy the following conditions: $ \min\{|S \cap V_1|, |S \cap V_2|\} = 1 $ when $ n $ is odd, and $ \min\{|S \cap V_1|, |S \cap V_2|\} \in \{1, 2, \frac{n}{2}\} $ when $ n $ is even.  Additionally, we identified infinitely many balanced bipartite graphs that meet these conditions. We demonstrated that there are also infinitely many balanced bipartite graphs $ G $ with $ \delta(G) = k $ and $ f(G) = \frac{|V(G)|}{2} + l $, where $ k \geq 2 $ and $ l \in \{1, 2\} $. So, to extend the results of this paper, we ask the following question: \\
\noindent 
{\bf Question.} What are the best conditions that can be established to identify all balanced bipartite graphs $ \mathcal{B} $ with $ f(\mathcal{B}) = n + 1 $?

It is a well-known fact that if $\mathcal{B}$ is a bipartite planar graph, then the number of edges $|E(\mathcal{B})|$ is at most $4n - 4$. Consequently, if $\mathcal{B}$ is a bipartite planar graph and $\delta(\mathcal{B})\geq\frac{n}{2} + 1$, it follows that $n \leq 4$. 
Furthermore, the inequality $\frac{10n}{8} > n + 1$ holds if and only if $n > 4$. Therefore, Theorem \ref{th1} cannot provide any counterexamples to the well-known conjecture regarding bipartite planar graphs mentioned in the introduction section. 
By a similar method, we can also see that the balanced bipartite graphs in Proposition \ref{pro1} and Theorems \ref{thh1} and \ref{th5} cannot provide any counterexamples for that conjecture either.

Finally, it is worth noting that based on the relationship $ f(\mathcal{B}) + \nabla(\mathcal{B}) = 2n $, the decycling version of the results presented in this paper can be directly derived.

\section*{Declaration of competing interest}
The authors declare that they have no known competing financial interests or personal relationships that could have
appeared to influence the work reported in this paper.
\section*{Acknowledgments}
We would like to thank the anonymous referees for their valuable and helpful comments and suggestions. The research was supported by the NSFC No.\ 12131013 and 12161141006.
\section*{Data availability}
No data was used for the research described in the article

\end{document}